\documentstyle[12pt,leqno,amssymb,graphics]{amsart}

\newtheorem{thm}{Theorem}[section]
\newtheorem{lemma}[thm]{Lemma}
\newtheorem{prop}[thm]{Proposition}

\theoremstyle{definition}

\theoremstyle{remark}

\begin{document}

\newcommand{\ct}{\cite}
\newcommand{\pr}{\protect\ref}
\newcommand{\su}{\subseteq}
\newcommand{\pa}{{\partial}}
\newcommand{\im}{{Imm(F,M)}}
\newcommand{\lm}{{\lambda}}
\newcommand{\hf}{{1 \over 2}}
\newcommand{\p}{{\pi_a}}
\newcommand{\ep}{{\epsilon}}
\newcommand{\tc}{{\mathfrak{T}}}
\newcommand{\Q}{{\Bbb Q}}
\newcommand{\R}{{\Bbb R}}
\newcommand{\Z}{{\Bbb Z}}
\newcommand{\E}{{{\Bbb R}^3}}
\newcommand{\B}{{\Bbb B}}
\newcommand{\G}{{\Bbb G}}

\newcommand{\A}{{\mathcal{A}}}
\newcommand{\C}{{\mathcal C}}
\newcommand{\I}{{\mathrm{Id}}}

\newcommand{\hM}{{\widehat{M}}}
\newcommand{\gr}{{\Z[\p]}}

\newcounter{numb}

\title{Immersions of surfaces into aspherical 3-manifolds}
\author{Tahl Nowik}
\address{Department of Mathematics, Bar-Ilan University, 
Ramat-Gan 52900, Israel.}
\email{tahl@@math.biu.ac.il}
\urladdr{http://www.math.biu.ac.il/$\sim$tahl}
\date{October 26, 2005}
\thanks{Partially supported by the Minerva Foundation}

\begin{abstract}
We study finite order invariants of null-homotopic immersions of a closed orientable
surface into an aspherical orientable 3-manifold. 
We give the foundational constructions, and classify all order one invariants.
\end{abstract}

\maketitle

\section{Introduction}
Finite order invariants of immersions of a closed orientable surface 
into $\E$ have been introduced in \ct{o},
where all order 1 invariants have been classified. Explicit formulae for the majority of order
1 invariants have been given in \ct{a},\ct{f}. All higher order invariants have been classified
in \ct{h}, and the analogue of all the above for non-orientable surfaces has appeared in \ct{n}.
A first step in the study of finite order invariants of immersions of surfaces into general 3-manifolds
appears in \ct{q} where one specific order 1 invariant has been studied. 

In the present work we study finite order invariants of null-homotopic immersions of a closed orientable
surface $F$ into an aspherical orientable 3-manifold $M$. 
The immersions being null-homotopic
will enable us to lift the immersions to the universal covering $\hM$ of $M$, and 
$M$ being aspherical will guarantee a simple form for $H_2(\hM - A)$ for discrete subset $A \su \hM$.
We develop the foundations for the study of finite order invariants
in our setting, as has been done in \ct{o} for the case $M=\E$ (and analogous to chord diagrams 
and the 1-term and 4-term relations in knot theory). We then classify all order 1 invariants. 

Revisions to this paper will appear at
www.math.biu.ac.il/$\sim$tahl/publications.html

\section{The universal covering space}\label{lb}

Let $M$ be an aspherical orientable 3-manifold, i.e. $\pi_n(M)=0$ for all $n \geq 2$
(a sufficient condition being that $\pi_2(M)=\pi_3(M)=0$).
Let $r:\hM \to M$ be the universal covering, so $\hM$ is contractible.
For $a \in M$ denote $\p=\pi_1(M,a)$. 
If $\gamma$ is a path in $M$ with initial point $a$ and if $x \in r^{-1}(a)$
then we denote the lift of $\gamma$ to $\hM$ with initial point $x$ by $\gamma^x$.
If $p \in M$ is another point and $\gamma$ is a path from $a$ to $p$ then we define a 
bijection $F^x_\gamma : \p \to r^{-1}(p)$ by
$F^x_\gamma(\phi) = (\phi * \gamma)^x (1)$, where $*$ denotes concatenation from left to right. 

We choose $a \in M$, $x \in r^{-1}(a)$ and for each $p \in M$ we choose a path
$\gamma_p$ from $a$ to $p$. We then label the points in each $r^{-1}(p)$ by elements
of $\p$ via the bijection $F^x_{\gamma_p}$.
We now check in what way this labeling depends on our choices.
So let $b \in M$, $y \in r^{-1}(b)$ and for each $p \in M$, $\delta_p$ a path from $b$ to $p$,
be another such choice.
In order to identify the elements of $\p$ with those of $\pi_b$ we need to make one additional
choice, the dependence on which will be apparent as well, namely we choose a path $\mu$ from $a$
to $b$ and identify $\p$ with $\pi_b$ by $\phi \mapsto \bar{\mu} * \phi * \mu$ where $\bar{\mu}$
is the path inverse to $\mu$ (i.e. $\bar{\mu}(t)=\mu(1-t)$).
Finally, let $s$ be the unique (up to homotopy) path in $M$ from $a$ to $b$ such that $s^x(1)=y$.
Now for each $p \in M$: 
$F^y_{\delta_p} (\bar{\mu} * \phi * \mu) = (\bar{\mu} * \phi * \mu * \delta_p)^y (1)
= (\bar{\mu} * \phi * \mu * \delta_p * \bar{\gamma_p} * \gamma_p)^{s^x(1)} (1)
= (s * \bar{\mu} * \phi * \mu * \delta_p * \bar{\gamma_p} * \gamma_p)^x (1)
= F^x_{\gamma_p} ( (s * \bar{\mu}) * \phi * (\mu * \delta_p * \bar{\gamma_p}))$.
So given the identification between $\p$ and $\pi_b$ determined by $\mu$, the labeling
of the points in $r^{-1}(p)$ have changed by left 
multiplication by $s * \bar{\mu}$ and right multiplication by 
$\mu * \delta_p * \bar{\gamma_p}$. Note that $s * \bar{\mu}$ depends only
on $a,b,x,y,\mu$ and so we get left multiplication by the \emph{same} element for all points $p$.
We summarize this in the following:

\begin{prop}\label{dep}
The labeling of the points of $r^{-1}(p)$ by elements of $\p$ is well defined
up to left multiplication by one common element in $\p$ for all $p$,
and right multiplication by an element in $\p$ which may depend on $p$.
\end{prop}

Once the dependence on choices has been established, we fix $a,x, \{\gamma_p\}_{p \in M}$ once and for
all, and so the labeling of the points in $r^{-1}(p)$ for each $p \in M$ is from now on fixed. 
Note that this labeling is necessarily \emph{not} locally constant.

There will be two families of bijections that we will have occasion to use, and we will now
see how they are expressed in terms of the labeling.
First let $D : \hM \to \hM$ be a deck transformation, then the restriction
of $D$ to each $r^{-1}(p)$ is a bijection onto itself.
Let $\psi \in \p$ be the unique element such that $\psi^x(1)=D(x)$.
Then for any $\phi \in \p$ and any $p \in M$, 
$D(F^x_{\gamma_p}(\phi)) = D((\phi * {\gamma_p})^x (1)) = (\phi * {\gamma_p})^{\psi^x(1)} (1) 
=(\psi * \phi * {\gamma_p})^x(1) = F^x_{\gamma_p}(\psi * \phi)$, that is, in terms of the labeling,
$D$ is given by common left multiplication by $\psi$. (By \emph{common} we mean as above,
that it is the \emph{same} $\psi$ for all $p$.)
Secondly, given $p,q \in M$ let $\delta$ be a path from $p$ to $q$, then $\delta$ 
defines a bijection $G_\delta : r^{-1}(p) \to r^{-1}(q)$ given by
$G_\delta(y) = \delta^y(1)$. We have for any $\phi \in \p$, 
$G_\delta (F^x_{\gamma_p}(\phi)) = G_\delta((\phi * \gamma_p)^x(1)) = 
\delta^{(\phi * \gamma_p)^x(1)}(1) = (\phi * \gamma_p * \delta)^x(1) = 
(\phi * (\gamma_p * \delta * \bar{\gamma_q}) * \gamma_q)^x(1) =
F^x_{\gamma_q}( \phi * (\gamma_p * \delta * \bar{\gamma_q}))$, that is, in terms of the labeling,
$G_\delta$ is given by right multiplication by $\gamma_p * \delta * \bar{\gamma_q}$.

It is also clear from the above calculations that any common left multiplying element, 
and any right multiplying element, can be realized by appropriately choosing $D$ and $\delta$ respectively.
We summarize this in the following:

\begin{prop}\label{dk}
\begin{enumerate}
\item The effect of a deck transformation is a common left multiplication, and any such common left
multiplication can be realized by a deck transformation.
\item The effect of $G_\delta$ is right multiplication, and any such right multiplication 
can be realized by some $\delta$ from $p$ to $q$.
\end{enumerate}
\end{prop}

We have noted that our labeling is not locally constant. By Proposition \pr{dk}(2) we see that any such 
discontinuity is always given by right multiplication.

Denote by $\gr$ the group ring of $\p$ with coefficients in $\Z$. We will in fact not use the full ring structure, 
but only the left and right actions of $\p$ on $\gr$. For each $p \in M$, 
$H_2(\hM - r^{-1}(p))$ is a free Abelian group with basis in one-to-one correspondence with 
$r^{-1}(p)$. This is true by a Mayer-Vietoris sequence, since $\hM$ is contractible. 
So via our labeling, for each $p \in M$ we may identify $H_2(\hM - r^{-1}(p))$ with $\gr$. 
Also note $H_2(\hM - r^{-1}(p)) = \pi_2(\hM - r^{-1}(p))$.

Any deck transformation induces an automorphism of $H_2(\hM - r^{-1}(p))$ for each $p$. It follows from
Proposition \pr{dk}(1) that in terms of the identification of $H_2(\hM - r^{-1}(p))$ with $\gr$, 
this automorphism (of $\gr$ as an Abelian group) 
is given by left multiplication by some common $\psi \in \p$. 
Similarly, given $p,q \in M$, a path $\delta$ from $p$ to $q$ naturally induces an isomorphism from
$H_2(\hM - r^{-1}(p))$ to $H_2(\hM - r^{-1}(q))$, and it follows from Proposition \pr{dk}(2) that
in terms of the identifications with $\gr$, this isomorphism is given by right multiplication by some element
of $\p$.

\section{Co-oriented AB equivalences}\label{imm}

Let $M$ be an oriented aspherical 3-manifold,
let $F$ be a closed oriented surface, and let $\A \su \im$ be a regular homotopy class
of immersions $i:F \to M$ which are null-homotopic. 
A CE point of an immersion $i:F \to M$ is a point of self intersection
of $i$ for which the local stratum in $\im$ corresponding to the 
self intersection, has codimension one. 
We distinguish four types of CEs which we name
$E, H, T, Q$. 
In the notation of \ct{hk} they are respectively
$A_0^2|A_1^+$, $A_0^2|A_1^-$, $A_0^3|A_1$, $A_0^4$.
The four types may be demonstrated by the following local models, where letting $\lm$ vary, 
we obtain a 1-parameter family of immersions which is transverse to the 
given codim 1 stratum, intersecting it at $\lm=0$.

$E$: \ \ $z=0$, \ \ $z=x^2+y^2+\lm$. 

$H$: \ \ $z=0$, \ \ $z=x^2-y^2+\lm$.

$T$: \ \ $z=0$, \ \ $y=0$, \ \ $z=y+x^2+\lm$. 

$Q$: \ \ $z=0$, \ \ $y=0$, \ \ $x=0$, \ \ $z=x+y+\lm$. 

See Figure \pr{hoi}, which corresponds to some small $\lm > 0$. If $R$ is one of the above four CE types,
then we denote by $|R|$ the number of sheets involved in the given configuration, that is,
$|E|=|H|=2$, $|T|=3$, $|Q|=4$.
Let $I_n \su \A$ be
the space of all immersions with precisely $n$ CEs, in particular, $I_0$ is the space of all stable immersions.

\begin{figure}[t]
\scalebox{0.7}{\includegraphics{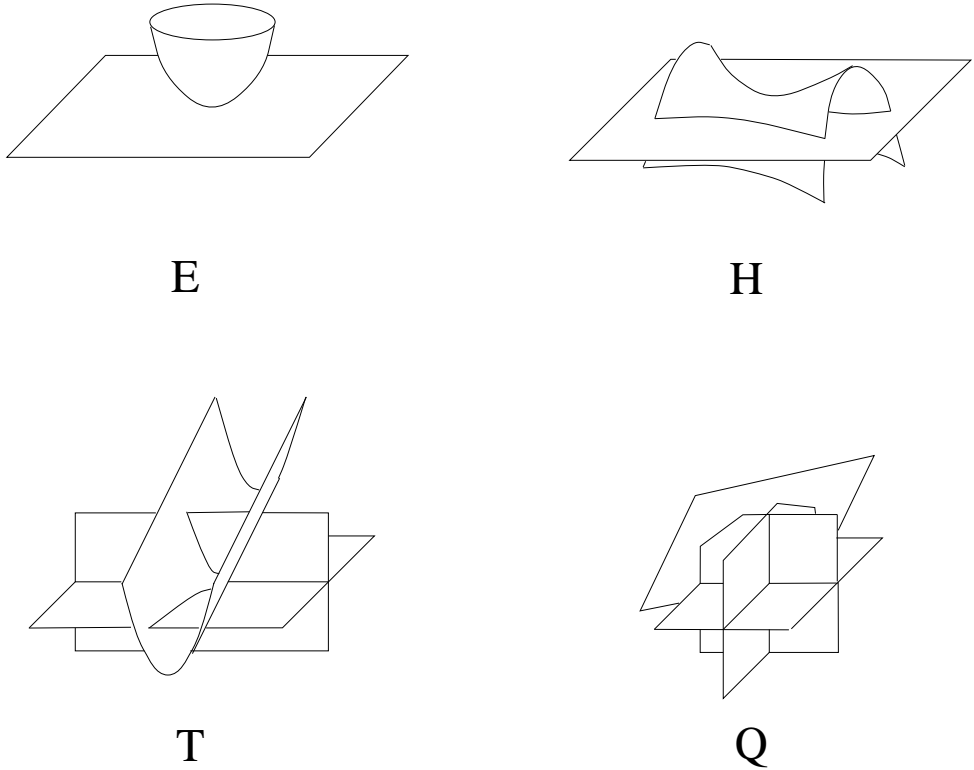}}
\caption{The four CE types}\label{hoi}
\end{figure}

A choice of one of the two sides of the local codim 1 stratum at a given point of the 
stratum, is represented by the choice of $\lm<0$ or $\lm>0$ in the formulae above. We will 
refer to such a choice as a 
\emph{co-orientation} for the configuration of the self intersection. 
For types $E$ and $T$, the configuration of the self intersection at the two sides 
of the stratum is distinct, namely, for $\lm<0$ there is an additional 
2-sphere in the image of the immersion, and we permanently
choose this side ($\lm<0$) as our positive side for the co-orientation. 
For types $H$ and $Q$, the configuration of the self intersection on the two 
sides of the strata is the same.

We will now further discuss the symmetries of our four configurations.
A symmetry of a CE configuration is an orientation preserving diffeomorphism 
from a neighborhood of the CE in $M$ to itself, which maps sheets onto sheets. 
A symmetry $h$ will be said to \emph{preserve} the co-orientation of the CE if the image under $h$ of motion 
into some side of the stratum, is again motion 
into the same side, (after identifying each sheet with its image).
The comment of the previous paragraph may now be stated as follows:
All symmetries of $E,T$ configurations preserve co-orientation, whereas configurations $H,Q$
admit co-orientation reversing symmetries.
We further notice, that configurations $E,H,T$ admit co-orientation preserving symmetries
that realize any given permutation of the sheets, whereas for configuration $Q$, any co-orientation
preserving symmetry induces an even permutation of the sheets, and any co-orientation reversing 
symmetry induces an odd permutation of the sheets. These facts may be seen by slightly resolving the
CE into one of the sides, and attempting to construct the diffeomorphism $h$ for the resolved configurations.

Given an immersion $i\in I_n$, an \emph{$n$-co-orientation} for $i$ is a choice 
of co-orientation at each of the $n$ CE points $p_1, \dots , p_n$ of $i$.
(This is called a temporary co-orientation in \ct{o},\ct{h},\ct{n}.)
A \emph{proper} $n$-co-orientation for $i \in I_n$ is an $n$-co-orientation where at each CE of type $E$ and $T$,
the co-orientation chosen is the permanent one mentioned above. (This definition is slightly weaker than that
in \ct{h}.) For immersion $i \in I_n$, and proper $n$-co-orientation $\tc$ for $i$,
let $p \in M$ be one of its $n$ CE points and we define $C_p(i,\tc)$ as the symbol
$R^{g_1\ep_1,\dots,g_{|R|}\ep_{|R|}}_d$ which is composed of:

(1) A symbol $R \in \{ E,H,T,Q \}$ which is the configuration of the given CE.

(2) A sequence of elements $g_1,\dots,g_{|R|} \in \p$ which is determined as follows:
Let $\hat{i}:F \to \hM$ be a lift of $i$. There are $|R|$ sheets of $\hat{i}(F)$ passing through points of
$r^{-1}(p)$ 
which in turn are labeled by elements $g_1,\dots,g_{|R|} \in \p$ (Section \pr{lb}).
If the CE is of type $E$,$H$ or $T$ then the order of the $g_j$s is arbitrarily chosen. 
As to $Q$ configuration, we choose once and for all a rule by which an orientation on a 
simplex determines an ordering on its four faces, up to an even permutation. Now, we resolve
the CE of type $Q$ positively with respect to $\tc$, 
i.e. slightly deform $i$ (not $\hat{i}$) near $p$,
to obtain a stable immersion in a neighborhood of $p$,
and write the elements $g_1,g_2,g_3,g_4$
in order determined by the orientation of $M$ restricted to the small simplex created by this resolution. 
(This ordering convention is related
to the symmetry properties discussed above).

(3) To each $g_j, 1 \leq j \leq |R|$, 
there corresponds a sign $\ep_j \in \{ +, - \}$ which is determined as follows:
We look at $i: F \to M$ (not the lift), and resolve the CE at $p$ according to $\tc$. 
For types $E,T,Q$ this creates a little sphere in the image of $i$, which bounds
a little 3-cell $V$ in $M$ (for $Q$ this is a simplex and for $E,T$ 
recall that $\tc$ is proper). 
Each $g_j$ we have, corresponds to a sheet involved in the CE. 
This sheet may have $V$ on its non-preferred side, 
(determined by the orientations of $F$ and $M$), 
in which case we set $\ep_j=+$, or otherwise
$V$ is on its preferred side, in which case $\ep_j=-$.
For type $H$ the region $V$ 
is not bounded by the local configuration,
but may still be defined, e.g. for $\lm>0$ in the formula for $H$ above,
$V$ will be a region consisting of points close to the origin and satisfying 
$0\leq z\leq x^2-y^2+\lm$. Now the signs $\ep_j$ are determined by the \emph{opposite} convention
than above, namely, 
a sheet having $V$ on its \emph{preferred} side will have $\ep_j=+$, and if
$V$ is on its \emph{non-preferred} side we set $\ep_j=-$.
The reason for the opposite convention is to obtain slightly nicer formulae in the end.
(It is consistent with the conventions in \ct{o}, chosen there in a similar way for the same reason).

(4) The subscript $d \in \gr$ which will be called the \emph{degree} of the CE,
is determined as follows: 
Let $\hat{i}:F \to \hM$ be the same lift as used in (2).
Wherever $\hat{i}(F)$ passes through a point of 
$r^{-1}(p)$, we push it slightly into the preferred side of $\hat{i}(F)$ in $\hM$, 
obtaining a map $F \to \hM - r^{-1}(p)$.
This map represents an element of $H_2(\hM - r^{-1}(p))$ which we identify with an element $d \in\gr$
(Section \pr{lb}).

If $p_1,\dots,p_n$ are the $n$ CEs of $i \in I_n$ then we define $C'(i,\tc)$
to be the $n$-tuple $(C_{p_1}(i,\tc),\dots,C_{p_n}(i,\tc))$, so 
$C'(i,\tc) \in \C'_n$ where $\C'_n$ is the set of all $n$-tuples of symbols of the form
$R^{g_1\ep_1,\dots,g_{|R|}\ep_{|R|}}_d$. 
In addition to the permanent choices of Section \pr{lb}, $C'$ also 
depends on the choice of lift $\hat{i}$, and various ordering choices. 
We thus define an equivalence relation on $\C'_n$ to be the equivalence relation generated by the 
following operations:

\begin{enumerate}
\item Any permutation of the $n$ symbols.
\item For $R=E,H,T$, any permutation of the $|R|$ elements 
$g_1\ep_1,\dots,g_{|R|}\ep_{|R|}$ (each pair $g_j\ep_j$ goes together).
\item For $R=Q$, any \emph{even} permutation of the four elements $g_1\ep_1,\dots,g_4\ep_4$.
\item Given $h \in \p$ replace each one of the $n$ symbols $R^{g_1\ep_1,\dots,g_{|R|}\ep_{|R|}}_d$ of the
$n$-tuple, by $R^{(hg_1)\ep_1,\dots,(hg_{|R|})\ep_{|R|}}_{hd}$ (the same $h$ for all $n$ symbols).
\item Given $h \in \p$, replace \emph{one} of the $n$ symbols $R^{g_1\ep_1,\dots,g_{|R|}\ep_{|R|}}_d$
by $R^{(g_1h)\ep_1,\dots,(g_{|R|}h)\ep_{|R|}}_{dh}$. 
(it is of course the same $h$ multiplying $g_j$ and $d$ on the right 
within one symbol, $h$ being allowed to vary only between the $n$ symbols).
\end{enumerate}

Denote by $\C_n$ the set of equivalence classes of $\C'_n$ under this equivalence relation.
Note that because of the \emph{common} left action in item (4), $\C_n$ is \emph{not} simply 
the set of unordered $n$-tuples of elements of $\C_1$.
We now define $C(i,\tc) \in \C_n$ to be the equivalence class of $C'(i,\tc)$.
So $C$ is well defined, independent of all choices, including the choice of lift $\hat{i}$
(Proposition \pr{dk}(1)). Denote by $I^P_n$ the set of all pairs $i,\tc$ where $i \in I_n$ and
$\tc$ is a proper $n$-co-orientation for $i$. Then $C : I^P_n \to \C_n$ and we claim:

\begin{lemma}\label{srj}
The map $C : I^P_n \to \C_n$ is surjective.
\end{lemma}

\begin{pf}
We are given an $n$-tuple $x$ of symbols of the form $R^{g_1\ep_1,\dots,g_{|R|}\ep_{|R|}}_d$.
Begin with any immersion $i \in \A$ and lift it to $\hat{i}:F \to \hM$. 
Choose points $p_1,\dots,p_n \in M$ as locations for the $n$ CEs we will construct. 
Deform $\hat{i}$  by regular homotopy
so that pieces of $F$ will pass the right points of $\hM$ above each $p_k$, to produce the correct labelings
appearing as superscripts in the symbols in $x$. Continue the deformation near each such point so that 
for the projection back to $M$ the right configurations will be created with the right signs, and inducing the
right ordering in case of CE of type $Q$.
Finally, move some other pieces of $F$ across the right points of 
$\bigcup_{1 \leq k \leq n} r^{-1}(p_k)$ to obtain the correct degree $d \in \gr$.
The projection back to $M$ gives the desired immersion.
\end{pf}

We also note that for types $H,Q$, 
the relation between the symbols obtained for a given co-orientation
of a CE, and that obtained for the opposite co-orientation, is as follows:
For $H$, all remains the same except for the signs $\ep_j$ which are all reversed.
For $Q$, all signs are reversed, and in addition an odd permutation is performed on 
$g_1\ep_1,\dots,g_4\ep_4$. 
(The degree $d$ remains unchanged since its definition does not involve the co-orientation.)
We call the symbol obtained from a symbol of type $H$ or $Q$
in this way, the \emph{reversed} symbol. 
(Recall that $C_p(i,\tc)$ is only defined for \emph{proper} $\tc$ and so 
for $E,T$ the co-orientation may not be reversed.)

We recall the definition of an AB equivalence, appearing in \ct{o}:
A regular homotopy between two immersions $i,j \in I_n$ is called an AB equivalence if 
it is alternatingly of type A and B, where 
\begin{enumerate}
\item $J_t:F\to M$ ($0\leq t \leq 1$)
is of type A if it is of the form
$J_t = U_t \circ i \circ V_t$ where 
$i:F\to M$ is an immersion and
$U_t: M \to M$, 
$V_t:F\to F$ are isotopies. 
\item $J_t:F\to M$ ($0\leq t \leq 1$)
is of type B if $J_0 \in I_n$ and there are  
little balls $B_1,\dots,B_n\su M$ centered at the $n$ CE points of $J_0$
such that $J_t$ fixes $U=(J_0)^{-1}(\bigcup_k B_k)$ 
and moves $F-U$ within $M - \bigcup_k B_k$.
\end{enumerate}

Given two immersions $i,j \in I_n$ and $n$-co-orientations $\tc,\tc'$ for $i,j$ respectively,
we now define \emph{co-oriented} AB equivalence, or CAB equivalence,
between $i,\tc$ and $j,\tc'$
to be an AB equivalence between
$i$ and $j$ which respects $\tc,\tc'$, i.e. if we carry $\tc$ at each CE of $i$
continuously along the AB equivalence, then we arrive at $j$ with $n$-co-orientation $\tc'$.

We now prove:

\begin{prop}\label{cab}
Let  $i,j\in I_n$ and $\tc,\tc'$ $n$-co-orientations for $i,j$ respectively. 
Then $i,\tc$ and $j,\tc'$ are CAB equivalent iff $C(i,\tc)=C(j,\tc')$.
\end{prop}

\begin{pf}
If $i,\tc$ and $j,\tc'$ are CAB equivalent, let $J_t$ ($0 \leq t \leq 1$)
be the given CAB equivalence, and $\tc_t$ ($0 \leq t \leq 1$), the $n$-co-orientation carried 
continuously along (so $\tc_1=\tc'$). If a lift $\hat{J}$ of $J$ 
is used for defining $C'(J_t,\tc_t)$ a each $t$, and all choices of orderings are also 
carried along continuously, then by Proposition 
\pr{dk}(2), $C'(J_t,\tc_t)$ may only change along the way by 
right multiplication. We may still have that $J_1$ is different 
than the lift of $j$ used to define $C'(j,\tc')$, 
in which case we will gain a common left multiplication.
The choices of ordering may also be different, but all together we get
$C(i,\tc)=C(j,\tc')$.

For the converse, 
assume $C(i,\tc)=C(j,\tc')$. Since our allowed ordering of superscripts, and equivalence of ordering, all correspond
to the possible symmetries of the configurations, there is an ambient isotopy $U_t :M \to M$ which brings the
CEs of $i$ onto the CEs of $j$, such that the $n$-co-orientation $\tc$ carried continuously along, 
coincides with $\tc'$, such that the orientations of the corresponding sheets match, (since the signs $\ep_j$ 
coincide,)
and such that all labelings in $\p$ and degrees in $\gr$ coincide 
up to right multiplication by an element in $\p$ and common left multiplication by an element in $\p$.
Since the orientations of corresponding sheets match, 
we may continue with a regular homotopy of the form $i \circ V_t$ for isotopy  $V_t : F \to F$,
until we have the same discs in $F$ participating in the CEs of $i$ and $j$, and the restriction of $i$
and $j$ to those discs coincides precisely, and we still have matching labeling and degrees up
to right and common left multiplication.
By Proposition \pr{dk} any right multiplication can be realized by dragging a CE around some loop in $M$ 
and back to its place to match $j$, and any common left multiplication can be realized by 
changing the choice of lift. So we perform such ambient isotopies and change of lift until we have that
all labelings and degrees coincide precisely. This means that now, not only $i,j$ but also their chosen lifts 
$\hat{i},\hat{j}$ coincide on the discs participating in the CEs, and that
if we slightly deform $i$ and $j$ by pushing each sheet of the CEs slightly into its preferred side, 
then $\hat{i},\hat{j} : F \to \hM - \bigcup_{1 \leq k \leq n} r^{-1}(p_k)$ represent the same element in
$H_2(\hM - \bigcup_{1 \leq k \leq n} r^{-1}(p_k))$.
(Note that 
$H_2(\hM - \bigcup_{1 \leq k \leq n} r^{-1}(p_k)) = \bigoplus_{1 \leq k \leq n}  H_2(\hM - r^{-1}(p_k))$,
where the projections are induced by inclusion.)
Form this point on we may proceed exactly as in the proof of \ct{o} Proposition 3.4, where instead of 
working in
$\E$ with set of designated points $\{p_1,\dots,p_n\}$, we work in $\hM$ with set of designated points
$\bigcup_{1 \leq k \leq n} r^{-1}(p_k)$. We then compose the obtained regular homotopy with $r$ obtaining the 
desired CAB equivalence in $M$.
\end{pf}

As already seen in \ct{o}, a CE of type $H$ or $Q$ 
may be CAB equivalent to itself with the opposite co-orientation, 
which means that the stratum corresponding to this CE in $\im$ is one sided. 
From Proposition \pr{cab} we see that this happens iff the corresponding symbol is equivalent
to its reversed symbol. For this to happen clearly the number of $+$ and $-$ signs should be equal, 
as this number is preserved under equivalence. So a one sided stratum may occur only for symbols of 
the form $H^{g_1+,g_2-}_d$ and $Q^{g_1+,g_2+,g_3-,g_4-}_d$. If all $g_j$ are equal then we clearly 
do get a symbol equivalent to its reversed, and so a one sided stratum. 
In \ct{o} where $\pi_1$ is trivial ($M=\E$), 
it is thus clear which are the one sided strata. But for general aspherical manifold this is intricately
related to the structure of $\pi_1(M)$. In the following two examples we show a case where the stratum 
for $H^{g_1+,g_2-}_d$ is one sided though $g_1 \neq g_2$, and then a case of an $M$ 
where the stratum is one sided only when $g_1=g_2$, though $\pi_1(M)$ is nontrivial.

Let $S$ be the connect sum of four projective planes, with $a,b,c,d \in S$ the standard four loops
for which $\pi_1(S)=\langle a,b,c,d \ | \ a^2b^2c^2d^2 = 1 \rangle$. Let $M$ be the 3-manifold
obtained from $S \times [0,1]$ by attaching a 2-handle to the loop $b d$ in $S \times \{ 1 \}$.
Then $\pi_1(M)=\langle a,b,c \ | \ a^2b^2c^2b^{-2}=1 \rangle$ and one can show that $M$ is 
indeed aspherical. We look at the symbol $H^{b^2 + , ab^2c -}_0 \in \C'_1$ (the subscript is the zero element of 
$\Z[\pi_1]$). We note that $b^2 \neq a b^2 c$, e.g. by looking at the quotient obtained 
by adding the relations $b=c=1$.   The reversed symbol 
is $H^{b^2 - , ab^2c +}_0 $. Left multiplication by $a$ and right multiplication by $c$ gives 
$H^{ab^2c - , a^2b^2c^2 +}_0 = H^{ab^2c - , b^2 +}_0$
(since $a^2b^2c^2=b^2$), which is equivalent back to $H^{b^2 + , ab^2c -}_0$ by permuting the 
superscripts. And so the $H^{b^2 + , ab^2c -}_0$ stratum is one sided.

On the other hand, for $M = S^1 \times \R^2$, if $x$ denotes a generator of
$\pi_1(M)$ and say $k < r$ then $H^{x^k+,x^r-}_0 \in \C'_1$ is not equivalent to the reversed 
symbol $H^{x^k-,x^r+}_0$ since the property of having the $+$ sign attached to the smaller
power of $x$ will be preserved under any left or right multiplication by element $x^s \in \pi_1(M)$.

\section{Finite order invariants}\label{fo}

Given an $n$-co-orientation $\tc$ for $i \in I_n$ and a subset $A\su \{p_1,\dots,p_n\}$,
we define $i_{\tc,A} \in I_0$ to be the immersion obtained from $i$ by resolving all CEs
of $i$ at points of $A$ into the 
positive side with respect to $\tc$,
and all CEs not in $A$ into the negative side.
Now let $\G$ be any Abelian group and let $f:I_0\to\G$ be an invariant, i.e. a function which is 
constant on each connected component of $I_0$.
Given $i\in I_n$ and an $n$-co-orientation $\tc$ for $i$,
$f^\tc(i)$ is defined as follows:
$$f^\tc(i)=\sum_{ A \su \{p_1,\dots,p_n\} } (-1)^{n-|A|} f(i_{\tc,A})$$
where $|A|$ is the number of elements in $A$.
If $\tc,\tc'$ are two $n$-co-orientations for the same immersion $i$ then
$f^\tc(i)=\pm f^{\tc'}(i)$ and so having $f^\tc(i)=0$ is independent of the $n$-co-orientation $\tc$.
An invariant $f:I_0\to\G$ is called \emph{of finite order} if 
there is an $n$ such that $f^\tc(i)=0$ for all $i\in I_{n+1}$.
The minimal such $n$ is called the \emph{order} of $f$.
The group of all invariants on $I_0$ of order at most $n$ is denoted $V_n = V_n(\G)$.

Let $f \in V_n$. If $i,j \in I_n$ and $\tc,\tc'$ are $n$-co-orientations for $i,j$ respectively, such that
$i,\tc$ and $j,\tc'$ are CAB equivalent, then $f^{\tc}(i) = f^{\tc'}(j)$, the argument being the same 
as in \ct{o} Proposition 3.8. And so we have by Lemma \pr{srj} and Proposition \pr{cab}, that 
each $f \in V_n$ induces a well defined function $\mu_n(f): \C_n \to \G$.
The map $f\mapsto \mu_n(f)$ induces an injection $\mu_n:V_n / V_{n-1} \to \C_n^*$ 
where $\C_n^*$ is the group of all functions from $\C_n$ to $\G$.

If $i \in I_n$ and $\tc,\tc'$ are two $n$-co-orientations for $i$ which differ at precisely
one CE, then for any $f \in V_n$, $f^\tc(i) = -f^{\tc'}(i)$.  If the CE where we have reversed the co-orientation 
is of type $H$ or $Q$, and located at $p \in M$, then we have already noticed the relation between
$C_p(i,\tc)$ and $C_p(i,\tc')$, and so we now get an equation that must be satisfied by any $g \in \C_n^*$ 
in order for it to lie in the image of $\mu_n$, namely, if $Z_2,\dots,Z_n$ are any $n-1$ symbols 
then for any symbol $H^{g_1\ep_1,g_2\ep_2}_d$, $g$ must satisfy
$g(H^{g_1\ep_1,g_2\ep_2}_d,Z_2,\dots,Z_n) = - g(H^{g_1\hat{\ep}_1,g_2\hat{\ep}_2}_d,Z_2,\dots,Z_n)$
where $\hat{\ep}_j$ denotes the sign opposite to $\ep_j$.
We will write such an equation in the short form 
$H^{g_1\ep_1,g_2\ep_2}_d = -H^{g_1\hat{\ep}_1,g_2\hat{\ep}_2}_d$. 
For $Q$ configuration we similarly get
$Q^{g_1\ep_1,g_2\ep_2,g_3\ep_3,g_4\ep_4}_d = -Q^{g_2\hat{\ep}_2,g_1\hat{\ep}_1,g_3\hat{\ep}_3,g_4\hat{\ep}_4}_d$
(note the transposition $1 \leftrightarrow 2$).
Also note that in the above two equations, the operation which produces the $n$-tuple 
of symbols on the right from the 
$n$-tuple of symbols on the left, is indeed well defined on the classes in $\C_n$.


Let $i \in \A$ be an immersion with a self intersection of local 
codim 2 at $p$ and 
$n-1$ additional self-intersections of local codim 1 (i.e. CEs) at $p_1,\dots,p_{n-1}$.
We look at a 2-parameter family of immersions which moves $F$ only in a neighborhood 
of $p$, such that the immersion $i$ corresponds to 
parameters $(0,0)$ and such that this 2-parameter family is transverse to the local 
codim 2 stratum at $i$.
In this 2-parameter family of immersions we look at a loop which 
encircles the point of 
intersection with the codim 2 strata,
i.e. a circle around the origin in the parameter plane.
This circle crosses the local codim 1 strata some $r$ times.
Between each two intersections we have an immersion in $I_{n-1}$ 
with the same $n-1$ CEs, at $p_1,\dots,p_{n-1}$.
At each intersection with the local codim 1 strata, an $n$th CE is added, 
obtaining an immersion in $I_n$. Let $i_1,\dots,i_r$ be the $r$ immersions in $I_n$ so 
obtained. For each $1 \leq k \leq r$, choose a proper $n$-co-orientation $\tc_k$ for
$i_k$, such that the co-orientation chosen for each $p_j$, $1 \leq j \leq n-1$ is the same in all
$\tc_1,\dots,\tc_r$. Let $e_k$, $k=1,\dots,r$ be $1$ or $-1$ according to whether we are passing 
the $n$th CE of $i_k$ in the direction of its co-orientation determined by $\tc_k$, or in the opposite direction,
respectively.
For an invariant $f$, it is easy to show that the following equation holds: 
$\sum_{k=1}^r e_k f^{\tc_k}(i_k) = 0$. 
Looking at $\mu_n:V_n / V_{n-1} \to \C_n^*$  we thus obtain additional equations that must be satisfied by 
a function in $\C_n^*$ in order for it to lie in the image of $\mu_n$.
We will now find all equations on $\C_n^*$ obtained in this way.
As above, the equations will be written in short form as equations on the symbols.
We may assume (by moving the codim 2 singular point $p$ if necessary), that there is a neighborhood 
$U$ of $p$ such that the labeling by elements of $\p$ is locally constant in $r^{-1}(U)$,
and such that all motion involved in our 2-parameter family of immersions takes place inside $U$.
And so we will be able to follow the labelings during our loop of immersions. 

The local codim 2 strata may be divided into six types which we name after the types of CEs 
appearing in a 2-parameter family of immersions, transverse to the given stratum:
$EH$, $TT$, $ET$, $HT$, $TQ$, $QQ$. In the notation of \ct{hk} they are 
respectively: $A_0^2|A_2$, $A_0^3|A_2$, $(A_0^2|A_1^+)(A_0)$, 
$(A_0^2|A_1^-)(A_0)$, $(A_0^3|A_1)(A_0)$, $A_0^5$.
Formula and sketch for local model for such strata,
the bifurcation diagrams, and the equations obtained, are as follows.
The symbol $\cdots$ represents a string of $g_j\ep_j$s of the appropriate length
(the same string for all appearances of $\cdots$ within the same diagram or equation).

\begin{figure}[h]
\scalebox{0.6}{\includegraphics{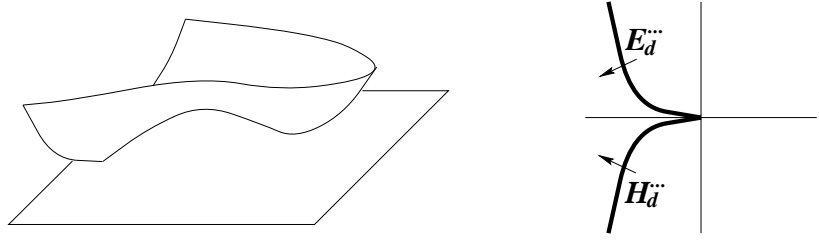}}
\caption{$EH$ configuration}\label{feh}
\end{figure}

$EH$: \ \ $z=0$, \ \ $z=y^2 + x^3+\lm_1 x + \lm_2$.
\begin{equation}\label{eeh}
0 = E^{\cdots}_d - H^{\cdots}_d
\end{equation}

\begin{figure}[h]
\scalebox{0.6}{\includegraphics{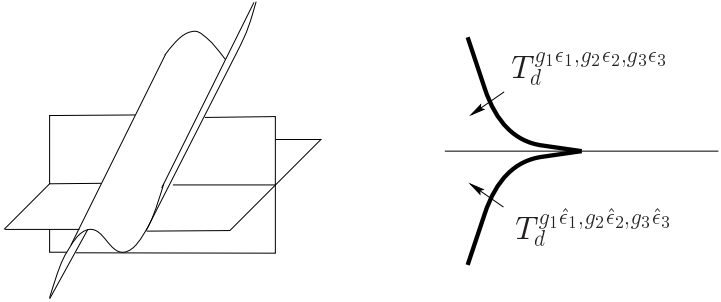}}
\caption{$TT$ configuration}\label{ftt}
\end{figure}

$TT$: \ \ $z=0$, \ \ $y=0$, \ \ $z=y+x^3+\lm_1 x  + \lm_2$.
\begin{equation}\label{ett}
0 = T^{g_1\ep_1,g_2\ep_2,g_3\ep_3}_d - T^{g_1\hat{\ep}_1,g_2\hat{\ep}_2,g_3\hat{\ep}_3}_d
\end{equation}

\begin{figure}[h]
\scalebox{0.6}{\includegraphics{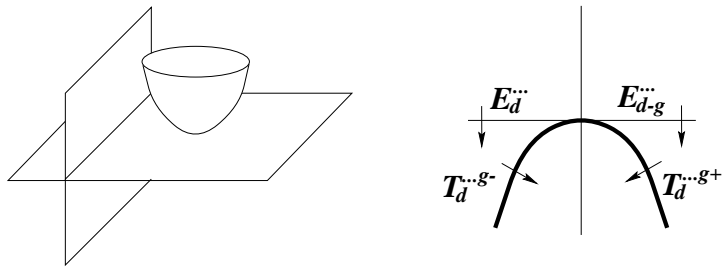}}
\caption{$ET$ configuration}\label{fet}
\end{figure}

$ET$: \ \ $z=0$, \ \ $x=0$, \ \ $z=(x-\lm_1)^2 + y^2 + \lm_2$.
\begin{equation}\label{eet}
0 = T^{\cdots g-}_d - T^{\cdots g+}_d - E^{\cdots}_{d-g} + E^{\cdots}_d
\end{equation}

\begin{figure}[h]
\scalebox{0.6}{\includegraphics{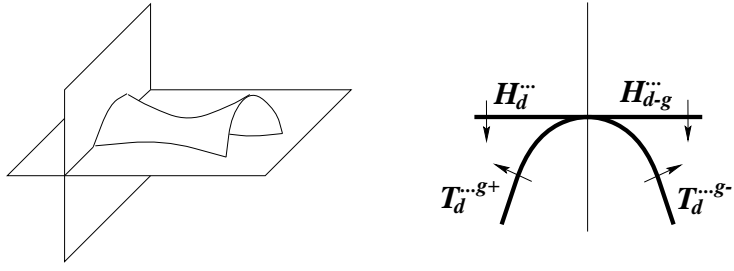}}
\caption{$HT$ configuration}\label{fht}
\end{figure}

$HT$: \ \ $z=0$, \ \ $x=0$, \ \ $z=(x-\lm_1)^2 - y^2 + \lm_2$.
\begin{equation}\label{eht}
0 = -T^{\cdots g+}_d + T^{\cdots g-}_d - H^{\cdots}_{d-g} + H^{\cdots}_d
\end{equation}

\begin{figure}[h]
\scalebox{0.6}{\includegraphics{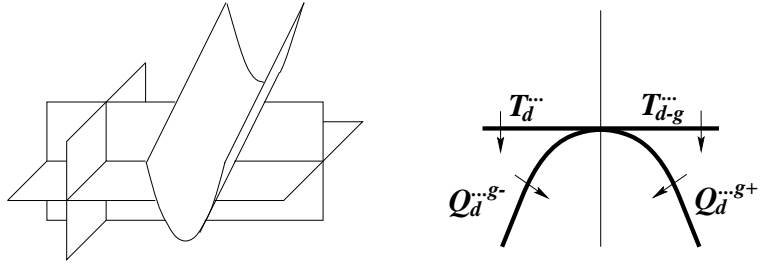}}
\caption{$TQ$ configuration}\label{ftq}
\end{figure}

$TQ$: \ \ $z=0$, \ \ $y=0$, \ \ $x=0$, \ \ $z=y+(x-\lm_1)^2 + \lm_2$.
\begin{equation}\label{etq}
0 = Q^{\cdots g-}_d - Q^{\cdots g+}_d - T^{\cdots}_{d-g} + T^{\cdots}_d
\end{equation}

$QQ$: This configuration is a quintuple point, i.e. five sheets passing through a
point, any three of which are in general position.
We model this by five planes passing through a point in $\E$.
A 2-parameter family may be constructed by fixing three of the planes, 
and allowing the other two planes to move parallel to themselves.
A loop in the shape of a rectangle around the origin of the parameter plane
may be constructed as follows: The three fixed planes $P_1,P_2,P_3$ intersect in a triple point 
$x$. Start with the plane $P_4$ very close to $x$, and $P_5$ somewhat further away. For the first
edge of the rectangle, move $P_4$ across
$x$ passing a first quadruple point. In this process a little simplex $S^-$ has vanished and a new little 
simplex $S^+$ has been created from the same four planes $P_1,\dots,P_4$.
For the second edge of the rectangle, 
move $P_5$ across the whole of $S^+$, passing four quadruple points on the way,
when crossing the four triple points of $S^+$.
For the third edge move $P_4$ back to its original place, thus crossing $x$ again, in the opposite direction,
recreating $S^-$ in place of $S^+$, this being the sixth quadruple point. Finally move 
$P_5$ back to place, across $S^-$, adding four 
more quadruple points, and so completing a total of ten quadruple points.
We co-orient each of these ten quadruple points according to the direction we are passing it along this loop,
and so each of the ten terms in the equation that we are producing, will have a $+$ sign in front of it.

We will call the first and third edge of the rectangle described above "short edges" (since they pass only
one quadruple point) and the second and fourth edges "long edges".
We will first find the relation between the symbols of two consecutive quadruple points along a
long edge, and then see that the same relation holds when passing to a short edge, and back.
So let $p,q$ be two vertices of the little simplex, which $P_5$ passes consecutively.
Let $e$ be the edge of the simplex connecting $p,q$. By affine transformation of $\E$ we may assume 
that $p = (0,0,-1), q=(0,0,1)$, that the two planes intersecting to create the edge $e$ are the planes 
$\{ x=0 \}$ and $\{ y=0 \}$, 
and that the plane $P_5$ which is in motion, crossing $p$ and then $q$, is the plane
$\{ z = t \}$ where $t$ increases in the range $-1-r < t < 1+r$ for some small $r$. 
For concreteness we look at the time $t=0$, and so $P_5 = \{ z=0 \}$.
As to the remaining two planes, one must pass $p$, 
and the other must pass $q$. 
They cannot be vertical planes (since we already have two vertical planes), 
and so they are given by $\{ z=-1+ax+by \}$ and $\{ z=1+cx+dy \}$. The coefficients
$a,b$ must both be nonzero since if say $b=0$ then we get three planes $\{ z=-1+ax \} $, $\{ z=0 \}$, 
$\{ x=0 \}$ which are not in general position (when parallelly moved to have a common point). Similarly $c,d$ 
are nonzero. 
By rotating the configuration around
the $z$ axis we may assume $a,b>0$. We claim that we must then get also 
$c,d >0$. Indeed, in order that $P_5$ will pass $p$ and $q$ consecutively, it must be that the other
two vertices of the simplex lie outside the region $-1 \leq z \leq 1$ of $\E$. 
These two points are $s_1 = \{ z=-1+ax+by \} \cap \{ z=1+cx+dy \} \cap \{ x=0 \}$
and $s_2 = \{ z=-1+ax+by \} \cap \{ z=1+cx+dy \} \cap \{ y=0 \}$. The $z$ component of $s_1$ is
$b+d \over b-d$, and this is outside the interval $[-1,1]$ iff $d>0$. 
(If $b=d$ then there is no intersection, which implies that the three planes are not in general position.)
In the same way $c>0$, using $s_2$.
And so the simplex $S^+_p$ created from the passage of $P_5$ across $p$ (bounded by the planes 
$\{ x=0 \}, \{ y=0 \}, \{ z=0 \}, \{ z=-1+ax+by \}$), lies in the octant $x \geq 0, y\geq 0, z \leq 0$,
whereas the simplex $S^-_q$ that is about to vanish when $P_5$ arrives at $q$ (bounded by the planes
$\{ x=0 \}, \{ y=0 \}, \{ z=0 \}, \{ z=1+cx+dy \}$), lies in the octant $x \leq 0, y\leq 0, z \geq 0$
which is the exact opposite octant. And so the sign $\ep_j$ with which each of the common planes
$\{ x=0 \}, \{ y=0 \}, \{ z=0 \}$ appears in $S^+_p$ and $S^-_q$ is the opposite sign. But since we
have co-oriented all quadruple points according to our direction of motion, the simplex $S^-_q$
is not the one with which we determine the signs for the symbol $Z_q$ of the quadruple point at $q$, but
rather the simplex $S^+_q$ which will be created \emph{after} we cross $q$.
As we have already noticed, the signs of the faces for $S^+_q$ are all opposite to the corresponding 
signs for $S^-_q$, and so finally the common planes with $S^+_p$ will have the same sign as in $S^+_p$. 

We apply a similar two step argument to determine the relation between the orderings of the faces in the 
two symbols $Z_p,Z_q$. We note that if some ordering of the faces of $S^+_p$ is consistent with the orientation
of $S^+_p$ (restricted from $\E$) then the ordering for the faces of $S^-_q$ obtained by simply replacing the
plane $\{ z=-1+ax+by \}$ with the plane $\{ z=1+cx+dy \}$ is \emph{in}consistent with the orientation of
$S^-_q$. And so this same ordering of these same four planes, now considered as the faces of $S^+_q$,
\emph{is} consistent with the orientation of $S^+_q$. 
And so for the superscripts of two consecutive symbols we finally get, that the $g_j\ep_j$ corresponding to
the plane which appears only in the first quadruple point and not the second (the plane $\{ z=-1+ax+by \}$
in our case) should be dropped, and exactly in its place should be written the $g_j\ep_j$ corresponding to
the new plane which participates in the second quadruple point but not in the first 
(the plane $\{ z=1+cx+dy \}$ in our case). The three $g_j\ep_j$s corresponding to the three common planes
(the planes $\{ x=0 \}, \{ y=0 \}, \{ z=\pm 1 \}$ in our case) must remain in their place, and with 
the signs $\ep_j$ unchanged.

We now show that when passing from the last quadruple point of a long edge to the quadruple point of the 
next short edge, or from that to the first quadruple point of the next long edge, then the same relation holds.
Indeed, looking say at the first and 
second quadruple point in our rectangle described above, let $y$ be the triple point
of $S^+$ which $P_5$ passes first along the long edge, for the occurrence of the second quadruple point. 
Let $P$ be one of the planes $P_1,\dots,P_4$ such that if $P$ is moved in the direction such that
$S^+$ increases in size, then the point $y$ will move toward $P_5$ and eventually cross it. So we can 
start from the same initial position as before, but instead of moving $P_4$ and then $P_5$, we move only
plane $P$, such that $S^-$ will vanish, passing the first quadruple point, $S^+$ is created and 
increases in size until $y$ crosses $P_5$ thus passing the second quadruple point.
So we have obtained the same two quadruple points, where now four of the planes are fixed and
only the plane $P$ is in motion. But this is precisely the motion analyzed above to obtain the 
relation between two consecutive quadruple points along a long edge.

Finally it remains to determine the degrees $d \in \gr$ of the ten symbols.
We do this again by analyzing the relation between the degrees of two consecutive 
quadruple points. Following the intersection point between $\{x=0\}, \{y=0\}, \{z=t\}$ as it continuously moves 
up from $p$ to $q$, we see that there is no change in the degree due to 
these three planes. 
(Recall that the labeling is locally constant above the region of interest.)
There may, on the other hand, be a change due to the planes $P_i = \{ z=-1+ax+by \}$ 
and $P_j = \{ z=1+cx+dy \}$. If the preferred side of $P_i$ is facing upward, which happens iff the corresponding 
sign $\ep_i$ appearing in the symbol $Z_p$ for the first quadruple point is $-$,
then when computing the degree $d_p$ for the first quadruple point, $P_i$ is slight pushed upward
so as to pass above the point $p$. But as $P_5$ moves upward to $q$, $P_i$ does not participate in the second 
quadruple point, and remains below the point $q$ which is where the degree is now computed. So in this case there 
is a change of $-g_i$ to the degree due to the plane $P_i$. If on the other hand the plane $P_i$ is facing 
downward, which happens iff the sign $\ep_i$ appearing in $Z_p$ is $+$, then the plane $P_i$ is 
pushed \emph{downward} when $d_p$ is computed, and again it will remain below the point $q$ when the degree 
of the second quadruple is computed, and so there will be no change in the degree due to $P_i$. The two cases 
may be written in one formula as follows. For sign $\ep \in \{ +,- \}$, define 
$$ | \ep | = 
\begin{cases}
1 & \text{if} \ \ep = + \\
0  &  \text{if} \ \ep = -
\end{cases}
$$
then the change in degree between the first and second quadruple point, due to the plane $P_i$ 
is $-|\hat{\ep_i}|g_i$. A similar analysis as to the effect of the plane $P_j$ on the degree will 
give $+|\ep_j|g_j$. Together we get that if $P_j$ is the plane missing from the first quadruple point, and $P_i$ 
appears there with sign $\ep_i$, and $P_i$ is missing from the second quadruple point, and $P_j$ appears there 
with sign $\ep_j$, then the change in the degree between the first and second quadruple point 
is $+|\ep_j|g_j-|\hat{\ep_i}|g_i$. 

Looking at the circle (or rectangle) of ten quadruple points, we see that the same 4-tuple of planes participates 
in a quadruple point twice, and with opposite co-orientation, and so the ordering of the corresponding $g_i\ep_i$ at 
the two times differ by an odd permutation. 
They occur at the precise opposite timing along the circle, i.e. five places apart, and so the ordering of 
quadruple points along the circle is that all five possible 4-tuples appear one after the other, and then they 
appear again in the same order, to complete the cycle of ten quadruple points. So given a quadruple point, either 
it or the matching quadruple point on the opposite side of the circle, have the property, that the ordering of 
superscripts may be chosen as
$g_1\ep_1,\dots,g_4\ep_4$ so that in the next quadruple point the plane corresponding to $g_1$ will be missing, in 
the following quadruple point the plane corresponding to $g_2$ will be missing, 
and so on. Finally, since the two matching quadruple points (i.e. those with the same 4-tuple of 
planes) occur with opposite co-orientation, the sign attached to each $g_i$ will be opposite 
(i.e. $g_i\ep_i$ is replaced by $g_i\hat{\ep}_i$). 
Since we have seen that if a plane is common to two consecutive quadruple points, it appears
in their superscripts with the same sign, and since a given plane will be missing exactly once between
its appearance in a given quadruple point, and the occurrence of the matching quadruple point, 
it follows that whenever a plane is missing from a quadruple point, then it appears 
in the following quadruple point with opposite sign than it does in the previous one.

To write our final formula, let us choose the first quadruple point and ordering $g_1\ep_1,\dots,g_4\ep_4$
as mentioned above, and we may write its degree as $d-|\ep_5|g_5$ where $g_5$ corresponds to the plane missing 
from the first quadruple point, and $\ep_5$ is the sign with which it will appear in the second quadruple point.
The first symbol is thus
$Q^{g_1\ep_1,g_2\ep_2,g_3\ep_3,g_4\ep_4}_{d-|\ep_5|g_5}$, and the combination of all analysis above produces the 
following equation:

\begin{multline}\label{eqq}
0 = Q^{g_1\ep_1,g_2\ep_2,g_3\ep_3,g_4\ep_4}_{d-|\ep_5|g_5}+ 
Q^{g_5\ep_5,g_2\ep_2,g_3\ep_3,g_4\ep_4}_{d-|\hat{\ep}_1|g_1}+ 
Q^{g_5\ep_5,g_1\hat{\ep}_1,g_3\ep_3,g_4\ep_4}_{d-|\hat{\ep}_2|g_2}+ 
Q^{g_5\ep_5,g_1\hat{\ep}_1,g_2\hat{\ep}_2,g_4\ep_4}_{d-|\hat{\ep}_3|g_3}+ 
Q^{g_5\ep_5,g_1\hat{\ep}_1,g_2\hat{\ep}_2,g_3\hat{\ep}_3}_{d-|\hat{\ep}_4|g_4} \\
+ Q^{g_4\hat{\ep}_4,g_1\hat{\ep}_1,g_2\hat{\ep}_2,g_3\hat{\ep}_3}_{d-|\hat{\ep}_5|g_5} + 
Q^{g_4\hat{\ep}_4,g_5\hat{\ep}_5,g_2\hat{\ep}_2,g_3\hat{\ep}_3}_{d-|\ep_1|g_1}+ 
Q^{g_4\hat{\ep}_4,g_5\hat{\ep}_5,g_1\ep_1,g_3\hat{\ep}_3}_{d-|\ep_2|g_2}+ 
Q^{g_4\hat{\ep}_4,g_5\hat{\ep}_5,g_1\ep_1,g_2\ep_2}_{d-|\ep_3|g_3}+ 
Q^{g_3\ep_3,g_5\hat{\ep}_5,g_1\ep_1,g_2\ep_2}_{d-|\ep_4|g_4}
\end{multline}
 
Now let $\Delta_n = \Delta_n(\G) \su \C_n^*$ be the subgroup consisting of all 
functions $g \in \C_n^*$ which satisfy the above equations (1) -- (6) together with our 
first two equations $H^{g_1\ep_1,g_2\ep_2}_d = -H^{g_1\hat{\ep}_1,g_2\hat{\ep}_2}_d$ and
$Q^{g_1\ep_1,g_2\ep_2,g_3\ep_3,g_4\ep_4}_d = -Q^{g_2\hat{\ep}_2,g_1\hat{\ep}_1,g_3\hat{\ep}_3,g_4\hat{\ep}_4}_d$.
Then we have obtained that the image of the injection $\mu_n : V_n / V_{n-1} \to \C_n^*$
is contained in $\Delta_n$. Finding the precise image of $\mu_n$ for all $n$ (as has been done for 
the case $M=\E$ in \ct{o},\ct{h}), would give a full classification of all finite order invariants. In the next section 
we will show that the image of $\mu_1$ is all $\Delta_1$, by this classifying all order one invariants.

\section{Order one invariants}\label{oo}

We define a ``universal'' Abelian group $\G_U$, by the Abelian group presentation
which takes as generators the elements of $\C_1$, and as relations, the set of all
equations used above to define $\Delta_n$ (and in particular $\Delta_1$). 
Note that for defining $\Delta_1$, an expression such as 
$H^{g_1\ep_1,g_2\ep_2}_d = -H^{g_1\hat{\ep}_1,g_2\hat{\ep}_2}_d$
was merely a short form for writing 
$g(H^{g_1\ep_1,g_2\ep_2}_d) = -g(H^{g_1\hat{\ep}_1,g_2\hat{\ep}_2}_d)$,
this being a condition on an element $g \in \C_1^*$ 
for being included in the subgroup $\Delta_1$, whereas, in the definition of $\G_U$,
$H^{g_1\ep_1,g_2\ep_2}_d = -H^{g_1\hat{\ep}_1,g_2\hat{\ep}_2}_d$
is an actual relation in the presentation of $\G_U$ by generators and relations.
We define the universal element $g^U \in \Delta_1(\G_U)$, to be the function that assigns to each
element in $\C_1$, the generator corresponding to it in $\G_U$. By definition of $\G_U$,
indeed $g^U \in \Delta_1(\G_U)$. We will now establish the existence of an order 1 invariant
$f^U : I_0 \to \G_U$ satisfying $\mu_1(f^U)=g^U$.
This will prove our desired result, that for any Abelian group $\G$, $\mu_1 : V_1(\G) \to \Delta_1(\G)$
is surjective. Indeed, given $g \in \Delta_1(\G)$, there is a (unique) homomorphism
$\varphi : \G_U \to \G$ such that $g = \varphi \circ g^U$, and we get 
$\mu_1 ( \varphi \circ f^U) = g$.

\begin{lemma}\label{B}
There exists $i_0 \in \A$ whose image is contained in a ball $B \su M$.
\end{lemma}

\begin{pf}
Let $i \in \A$ be some immersion and $H_t:F \to M$ a null-homotopy from $i$
to a constant map $k:F \to M$ with $k(F)=p \in B$. Using the parallelizibility
of $M$ we extend $H$ to a homotopy of bundle monomorphisms, from $di$ to a bundle monomorphism
$b$ which covers $k$. By the Smale-Hirsch Theorem 
applied to $Imm(F,B)$, there is an immersion $i_0:F \to B$ with $di_0$ homotopic to $b$ as bundle
monomorphisms. Now by the Smale-Hirsch Theorem applied to $Imm(F,M)$, $i$ is regularly homotopic
to $i_0$, and so $i_0 \in \A$.
\end{pf}

Fix an immersion $i_0 \in I_0$ with image contained in a small ball $B \su M$, as
provided by Lemma \pr{B}. 
For any $i \in I_0$ take a generic regular homotopy $J_t$ from $i_0$ to $i$,
and to each CE which occurs along $J_t$ choose a proper co-orientation.
Let $f^U(i)$ be the sum of the elements in $\G^U$ corresponding to the 
co-oriented CEs occurring along $J_t$, each taken with a sign $\pm$ according to whether
$J_t$ passes it in the direction of its co-orientation. An opposite choice of co-orientation 
for a CE along $J_t$ (relevant for types $H,Q$), will produce the negative element in $\G_U$, but 
this element will appear with opposite sign in the sum and so this sum is well defined, for
given $J_t$. 
We would like to show that in fact $f^U$ is
independent of the regular homotopy $J_t$, and so $f^U$ is a well defined invariant.
This is equivalent 
to showing that this sum is 0 along any closed $J_t$, i.e. $J_t$ such that $J_0=J_1=i_0$.
Once we know $f^U$ is well defined, it is clear that it is of order 1,
and $\mu_1(f^U)=g^U$.

The fact that the sum of values is 0 along any null-homotopic loop follows from
the definition of $\Delta_1$ and $\G_U$ as demonstrated in \ct{o}, 
and so the sum along loops induces a well defined homomorphism $\phi : \pi_1(\A,i_0) \to \G_U$.
We must show $\phi=0$.

Fix a disc $D \su F$, and let $S: \A \to Imm(D,M)$ be given by restriction,
$S(i)=i|_D$. Then $S$ is a fibration with fiber $\A_D = \{ i \in \A \ : \ i|_D = i_0|_D \}$.
So we get an exact sequence 
$$\pi_1(\A_D,i_0) \stackrel{inc_*}{\longrightarrow} \pi_1(\A , i_0) 
\stackrel{S_*}{\longrightarrow} \pi_1(Imm(D , M ) , i_0|_D ).$$
Let $K \su \pi_1(\A , i_0)$ be the subgroup consisting of loops obtained by composing $i_0 : F \to B$ 
with a motion of $B$ in $M$ through embeddings, beginning and ending with the inclusion of $B$ in $M$.  
It is clear that $S_*$ maps $K$ onto $\pi_1(Imm(D , M ) , i_0|_D )$, and so the elements of $K$ and 
those coming from $\pi_1(\A_D,i_0)$ generate $\pi_1(\A , i_0)$. Now $\phi$ is 0 on elements of $K$, since no 
CEs occur along such loops, and so it remains to show that $\phi$ is 0 on 
$\pi_1(\A_D,i_0)$. It is shown in the proof of \ct{q} Theorem 3.4, that $\pi_1(\A_D,i_0)$ is generated by
a loop which only moves a small disc $U \su F$ which is disjoint from $D$, 
and $U$ moves only within a small ball $B'$.
In our case we may take $B' \su B$, and so we finally need only to check the value of $\phi$ 
on loops moving $F$ in $B$. If we take our labeling to be locally constant above $B$,
then we can have all lifts used to define our symbols to be contained in the same lift of $B$,
constantly labeled by some fixed element $h \in \p$. 
And so all symbols will be of the form $R^{h\ep_1,h\ep_2,\dots}_{mh}$.
By a permutation of the superscript, which is even in case $R$ is $Q$,  we may assume that all
$+$ signs appear first, and so the superscript is characterized only by the \emph{number} of $+$ signs,
which is precisely the way the superscript of a symbol is defined in \ct{o} for immersions 
$F \to B = \E$.
The subscript is characterized by the coefficient $m$ of $mh$ which is also simply the degree
as defined in \ct{o} for immersions $F \to B$. And so $\phi$ being 0 on all loops moving
$F$ in $B$, follows from the fact that $f^U$ appearing in \ct{o} for the case $M=\E$, is well defined. 
Indeed, for $f^U$ of \ct{o}, the sum of symbols along a loop of immersions 
in $B$, is 0 in the group $G_U$ of \ct{o}. 
We have noticed that our symbols of the form $R^{h\ep_1,h\ep_2,\dots}_{mh}$
correspond to the symbols of \ct{o}, 
but furthermore, any equation on the symbols appearing in \ct{o}, 
has a corresponding equation on our symbols of the form 
$R^{h\ep_1,h\ep_2,\dots}_{mh}$, 
by considering the same codim 2 configuration in $B$. 
And so the sum in our present $G_U$ must also be 0. 
This completes the proof of: 

\begin{thm}
Let $F$ be a closed orientable surface, $M$ an aspherical orientable 3-manifold,
and $\A$ a regular homotopy class of null-homotopic immersions of $F$ into $M$.
Then for any Abelian group $\G$, the injection $\mu_1 : V_1 / V_0 \to \Delta_1$ is surjective.
\end{thm}

We conclude with the following remark. Our classification is in terms of the group $\Delta_1(\G)$
which as we have seen is the same as the group $Hom(\G_U , \G)$. In \ct{o} where $M=\E$ and so
$\pi_1(M)=0$, the structure of $\G_U$ is completely understood and accordingly the group
of order 1 invariants. In the general case presented here, the structure of
$\G_U$ depends on the structure of $\pi_1(M)$, most notably through the equivalence relation which
produces $\C_1$ from $\C'_1$. The complexity of this dependence may be seen in the 
discussion following Proposition \pr{cab}, where the question is whether a given symbol is equivalent
to its reversed symbol. And so a closed explicit classification as in \ct{o} 
(not to mention explicit formulae for the invariants, as in \ct{f}), will require extensive
analysis of $\pi_1(M)$ for any given $M$.

\end{document}